\documentclass[12pt]{article}
\usepackage{amsxtra,amssymb,amsthm,amsmath,latexsym}

\theoremstyle{plain}

\def\R{{\mathbb R}}

\def\oH{{\overset{\circ}{H}}}
\def\oH1{{\overset{\circ}{H}\kern-.02in{}^1}}

\def\bee{\begin{equation*}}
\def\eee{\end{equation*}}
\def\be{\begin{equation}}
\def\ee{\end{equation}}

\begin{document}

\title{ 
When does a double-layer potential equal to a single-layer one? 
}

\author{Alexander G. Ramm\\
 Department  of Mathematics, Kansas State University, \\
 Manhattan, KS 66506, USA\\
ramm@ksu.edu\\
http://www.math.ksu.edu/\,$\sim $\,ramm
}

\date{}
\maketitle\thispagestyle{empty}

\begin{abstract}
 Let $D$ be a bounded domain in $\R^3$ with a closed, smooth, connected boundary $S$, $N$ be the outer unit normal to $S$, $k>0$ be a constant,  $u_{N^{\pm}}$ are the limiting values of the normal derivative of $u$ on $S$ from $D$, respectively $D':=\R^3\setminus D$; $g(x,y)=\frac{e^{ik|x-y|}}{4\pi |x-y|}$, $w:=w(x,\mu):=\int_S g_{N}(x,s)\mu(s)ds$ be the double-layer
 potential, $u:=u(x,\sigma):=\int_S g(x,s)\sigma(s)ds$ be the single-layer potential.
 
 In this paper it is proved that for every $w$ there is a unique $u$, such
 that $w=u$ in $D$ and vice versa.
 
  Necessary and sufficient conditions are
 given for the existence of $u$ and the relation $w=u$ in $D'$, given $w$ in $D'$,  and for  the existence of $w$ and the relation $w=u$ in $D'$, given $u$ in $D'$. 
 
 \end{abstract}

\section{Introduction}\label{S:1}
 Let $D$ be a bounded domain in $\R^3$ with a closed, smooth, connected boundary $S$, $N=N_s$ be the outer unit normal to $S$ at the point $s\in S$, $k>0$ be a constant,  $g(x,y)=\frac{e^{ik|x-y|}}{4\pi |x-y|}$, 
 $w:=w(x,\mu):=\int_S g_{N}(x,s)\mu(s)ds$ be the double-layer potential,  $w^{\pm}$ are the limiting values of $w$ on $S$ from $D$, respectively, from $D'$,
 $u:=u(x,\sigma):=\int_S g(x,s)\sigma(s)ds$ be the single-layer potential,  $u_{N^{\pm}}$ are the limiting values of the normal derivative of $u$ on $S$ from $D$, respectively $D':=\R^3\setminus D$, the overbar, e.g. $\bar{\sigma}$, denotes the complex conjugate of $\sigma$,
 $ H^0:=L^2(S)$, $H^1:=W^1_2(S)$ is the Sobolev space, $\mu\in H^0$.  We write iff for if and only if and use the known formulas for the
 limiting values of the potentials on $S$, see \cite{R635}, pp. 148, 153:
 \be\label{e1}
 w^{\pm}=0.5(A'\mp I)\mu; \quad  u_N{^\pm}=0.5(A\pm I)\sigma;\quad w_{N}^{+}=w_{N}^{-}; \quad Q(\sigma):=\int_S g(t,s)\sigma(s)ds,
 \ee
where $A'\mu=\int_Sg_{N_t}(t,s)\mu(s)ds$, $A\sigma=\int_Sg_{N_s}(t,s)\sigma(s)ds$.
In this paper it is proved that for every $w$ in $D$ there is a unique $u$, such
 that $w=u$ in $D$, and vice versa.
 
  Necessary and sufficient conditions are
 given for $w=u$ in $D'$ and for $u=w$ in $D'$. 
 
 Let us state these results:
 
 {\bf Theorem 1.} {\em For every $w$, defined in $D$, there is a unique $u$ such that $w=u$ in $D$, and vice versa.
 
  For every $w$, defined in $D'$, there is a unique $u$ such that $w=u$ in $D'$ iff
\be\label{e2}
\int_Sw^{-}pds=0 \quad \forall p\in N(Q). 
\ee  
  For every $u$, defined in $D'$,  there is a unique $w$ such that $w=u$ in $D'$ iff
\be\label{e3}
\int_Su rds=0 \quad \forall r\in N(A+I). 
\ee  
}

In Section 2 proofs are given. In Section 3 it is proved that $Q: H^0\to H^1$
is a Fredholm operator.

\section{Proofs}\label{S:2}
a) Assume that $w=w(x,\mu)$ is given in $D$. Let us prove that $u=u(x,\sigma)$ exists such that $u=w$ in $D$, and $u$ is uniquely defined by $w$.
First, let us prove the last claim. Suppose $u_1=w$ and $u_2=w$ in $D$.
Let $u_1-u_2:=u$, $u=\int_Sg(x,s)\sigma  ds$ in $D$, $\sigma=\sigma_1-\sigma_2$. Then $u|_S=0$,
$(\nabla^2+k^2)u=0$ in $D'$ and $u$ satisfies the radiation condition,
so $u=0 $ in $D'$. Thus, $u=0$ in $D\cup D'$, $u_1=u_2$, and the claim is proved. 

Let us now prove the existence of $u$ such that $w=u$ in $D$. One has
\be\label{e4}
w^{+}=0.5(A'-I)\mu=u|_S=Q(\sigma) 
\ee  
This is an equation for $\sigma$, $w^{+}$ is given.  This equation is equivalent to
\be\label{e4'}
w^{+}=(I+Q_1Q_0^{-1}) Q_0(\sigma)=(I+Q_1Q_0^{-1})\eta, \quad \eta:=Q_0 (\sigma), \quad \sigma=Q_0^{-1}\eta.
\ee  
The operator $Q_1Q_0^{-1}$ is compact in $H^0$, see Section 3.
Therefore a necessary and sufficient condition for the solvability of equation \eqref{e4'}, and the equivalent equation \eqref{e4}, is:  
\be\label{e5}
\int_Sw^{+}  \bar{\eta} ds=0 \quad \forall  \eta \in N((I+Q_1Q_0^{-1})^\star),
\ee  
where $N(B)$ is the null space of the operator $B$ and
$B^\star$ is the adjoint operator to $B$ in $H^0$, $B=I+Q_1Q_0^{-1}$,
$B$ is of Freholm type in $H^0$. 
The kernel $g(t,s)$ of $Q$, the function 
$g(t,s)=\frac{e^{ik|t-s|}}{4\pi |t-s|}$, is symmetric: $g(t,s)=g(s,t)$.
Therefore, the kernel of $Q^\star$ is $\bar{g}(t,s)$.
 Clearly, $ (I+Q_1Q_0^{-1})^{\star} \eta =0$ iff  $ (I+\bar{Q}_1\bar{Q}_0^{-1}) \eta =0$, or, taking the complex conjugate,
 $ (I+Q_0^{-1}Q_1)\bar{\eta}=0$. Applying the operator $Q_0$ to this equation one gets an equivalent equation
 $(Q_0+Q_1)\bar{\eta}=0$ since $Q_0$ is an isomorphism.
 Let $u=u(x,\bar{\eta})$. Then $u|_S=0$.
Since $(\nabla^2+k^2)u=0$ in $D'$ and $u$ satisfies the radiation condition, it follows that $u=0$ in $D'$ and $\bar{\eta}=u^{+}_N-u^{-}_N=u^{+}_N$. Therefore, using the Green's formula, one obtains: 
$$\int_Sw^{+}  \bar{\eta} ds=\int_S w^{+}_N uds=0,$$
and, since $u=0$ on $S$, it follows that
\be\label{e6}
\int_Sw^{+}  \bar{\eta} ds=0 \quad \forall  \eta \in N((I+Q_1Q_0^{-1})^{\star}). 
\ee 
Thus, the necessary and sufficient condition  \eqref{e5} for the solvability of \eqref{e4} is always satisfied. Therefore,  $u(x,\eta)=w(x, \mu)$ in $D$,
$Q_0^{-1}\eta=\sigma$. 

b) Asume now that $u$ is given in $D$ and prove the existence of  a unique $w$ such that $u=w$ in $D$. 
Let us prove that  $w$ is uniquely determined by $u$ if  $u=w$ in $D$. 
Indeed, assume that there are two $w_j$, $j=1,2,$ such that  $u=w_j$ in $D$. Then $w:=w_1-w_2=0$ in $D$. Therefore $w_{N^+}=0$ on $S$.
It is known (see \cite{R635}, p. 154) that $w_{N^+}=w_{N^-}$, so $w_{N^-}=0$.
Therefore, $(\nabla^2+k^2)w=0$ in $D'$,  $w_{N^-}=0$ on $S$, and $w$
satisfies the radiation condition at infinity. This implies $w=0$ in $D'$, so
$w=0$ in $D\cup D'$. Therefore $\mu=w^--w^+=0$, so $w_1=w_2$ if
 $u=w_j$ in $D$.  We have proved that $w$ is uniquely determined by $u$
 if $u=w$ in $D$.

 Now the existence of the solution $\mu$ to
 equation \eqref{e4} 
and the relation $u=w$ in $D$ should be proved. The operator $A'-I$ is Fredholm in $H^0$, so a necessary and sufficient condition for the equation 
\eqref{e4} to be solvable  is:
\be\label{e7}
\int_S Q(\sigma) \bar{h} ds=0 \quad \forall  h \in N((A'-I)^{\star}). 
\ee  
One has $(A'-I)^{\star} h=(\bar{A}-I)h=0$ iff $(A-I)\bar{h}=0$.
If  $(A-I)\bar{h}=0$, then $u_N^{-}(s, \bar{h})=0$, so $u(s,\bar{h})=0$ in $D'$.  Note that $u(s,\bar{h})=Q(\bar{h})$. Therefore,
$Q(\bar{h})=0$ in $D'$. Since $Q$ is a symmetric operator in $H^0$, one has: $\int_SQ(\sigma)\bar{h} ds=\int_S \sigma Q(\bar{h}) ds=0$.  Consequently, condition \eqref{e7} is always satisfied, the solution $\mu$ to equation \eqref{e4} exists and  $u=w$ in $D$.

c)  Assume that $w$ is given in $D'$. Let us prove that $u$ exists such that 
$u=w$ in $D'$ iff 
\be\label{e8}
\int_Sw^{-}\bar{p} ds=0   \quad \forall  p \in N(Q^{\star}). 
\ee  
Consider the equation for $\sigma$:
\be\label{e9}
w^{-}=Q(\sigma).
\ee
 The operator $Q:H^0\to H^1$ is Fredholm-type. Thus, a necessary and sufficient condition for the solvability of the above equation is equation \eqref{e8}. If $\sigma$ solves \eqref{e9}, then $u(x,\sigma)=w$
 in $D'$ because the value of $u$ on $S$ determines uniquely $u$ in $D'$. 

d) Assume now that $u$ is given in $D'$.  Let us prove that $w$ exists such that  $u=w$ in $D'$ iff 
\be\label{e10}
\int_S u \bar{p} ds=0  \quad \forall  p \in N(\bar{A}+I). 
\ee  
Note that $ p \in N(A+I)$ iff $\bar{p}\in N(\bar{A}+I)$. The equation for $\mu$, given $u$ in $D'$, is:
\be\label{e11}  0.5 (A'+I)\mu=u.
\ee
Since the operator $A'+I$ is Fredholm in $H^0$, a necessary  and sufficient condition for the solvability of \eqref{e11} for  $\mu$ is:
\be\label{e12}
\int_S u \bar{p} ds =0 \quad \forall p\in N((A'+I)^{\star}).
\ee
If $p\in N((A'+I)^{\star})$, then $p\in N(\bar{A}+I)$, so condition \eqref{e10} is the same as \eqref{e12}.  As in section c), the relation $u=w$ in $D'$ is a consequence of the fact that $u=w^{-}$ on $S$.

Theorem 1 is proved. \hfill$\Box$

\section{Auxiliary lemma}\label{S:3}
Recall that $H^1$ is the Sobolev space on $S$, $H^0=L^2(S)$.

{\bf Lemma 1.} {\em  The operator $Q=Q_0+Q_1: H^0\to H^1$ is of Fredholm-type, where $Q_0$ is the operator
with the kernel $\frac 1 {4\pi |t-s|}$ and $Q_1$ has the kernel 
$\frac {e^{ik|t-s|}-1}{4\pi |s-t|}$. 
The operator $Q_0$ is an isomorhism of $H^0$ onto $H^1$, which has a continuous inverse. The operator $Q_1Q_0^{-1}$ is compact in $H^0$. }

{\em Proof.} Let us check that $Q_0: H^0 \to H^1$ is an isomorphism.  The Fourier transform of the kernel 
$\frac 1 {4\pi |x-y|}$ is positive: 
$\int_{\R^3}\frac {e^{i\xi \cdot x}}{4\pi |x|}dx=\frac 1 {|\xi|^2}$.
So, $Q_0: H^0 \to H^1$ is injective. 
The kernel of the operator $Q_1$ is smooth enough for $Q_1: H^0 \to H^1$  to be compact.
 Let us check that $Q_0: H^0 \to H^1$  is surjective. Let $f\in H^1$ and
 $Q_0\sigma=f$. Then $u:=u(x,\sigma)=\int_Sg(x,s)\sigma ds$ solves the 
 problem: $\nabla^2 u=0$ in $D$, $u|_S=f$. By the known elliptic estimates
(see, e.g., \cite{GT}) one has  $\|u\|_{H^{3/2}(D)}\le \|u\|_{H^1(S)}$.
Therefore, $\nabla u\in H^{1/2}(D)$ and, by the trace theorem, $u|_S\in H^0(S)$. This proves surjectivity of $Q_0: H^0 \to H^1$. Thus, $Q_0$ is an isomorphism of $H^0$ onto $H^1$ which has a continuous inverse. The sum of an isomorphism $Q_0$ and a compact operator $Q_1$ is a Fredholm operator, see, e.g., \cite{KG}. The operator $Q_1Q_0^{-1} $ is compact in $H^0$ because the kernel of $Q_1$ is sufficiently smooth. Although the 
  operator $Q_1Q_0^{-1}$ is defined on a dense subset $H^1$ of $H^0$,
  but since this operator is bounded in $H^0$ its closure is a bounded
  operator in $H^0$. Since the kernel of $Q_1$ is $O(|s-t|)$, the kernel of
  $Q_1Q_0^{-1}$ is a continuous function of $|s-t|$ and the surface $S$
is a compact set. Therefore, the operator  $Q_1Q_0^{-1} $ is compact in $H^0$.
 
Lemma 1 is proved. \hfill$\Box$

{\em Competing Interests:}  the author does not have such interests; this work 
was not supported financially by any source.

\end{document}